\documentclass[12pt]{article}
\usepackage{latexsym, amsmath, amssymb, amsfonts, amscd}
\usepackage{graphics}
\usepackage[dvips]{epsfig}

\newtheorem{theorem}{Theorem}
\newtheorem{proposition}{Proposition}
\newtheorem{lemma}{Lemma}

\newcommand{\call}{{\cal L}}
\newcommand{\calo}{{\cal O}}

\newcommand{\adL}{\mbox{\rm ad}_{\Lambda}}

\def\lcf{\lbrack\! \lbrack}
\def\rcf{\rbrack\! \rbrack}

\def\dbar{\overline\partial}

\newcommand{\CC}{\mathbb{C}}

\def\dbar{\overline\partial}

\def\oomega{\overline\omega}
\def\Oomega{\overline\Omega}

\newcommand{\lie}[1]{\mathfrak{#1}}
\newcommand{\lieg}{\mathfrak{g}}
\newcommand{\liet}{\mathfrak{t}}
\newcommand{\liec}{\mathfrak{c}}
\newcommand{\lieh}{\mathfrak{h}}

\newcommand{\lbra}[2]{\lcf #1, #2 \rcf}

\newcommand{\bproof}{\noindent{\it Proof: }}
\newcommand{\eproof}{\hfill \qed \vspace{0.2in}}
\def\qed{\rule{2.3mm}{2.3mm}}

\begin{document}
\title{\bf Holomorphic Poisson Cohomology}
\author{
Zhuo Chen\thanks{Address: Department of Mathematical Sciences, Tsinghua University, Beijing, P.R.C.} \ \
Daniele Grandini\thanks{ Address:
Department of Mathematics, University of New Mexico at Albuquerque,
    Albuquerque, NM 87131, U.S.A.} \ \
Yat-Sun  Poon\thanks{ Address:
    Department of Mathematics, University of California at Riverside,
    Riverside, CA 92521, U.S.A.. Email: ypoon@ucr.edu.} }
\maketitle
\begin{abstract} Holomorphic Poisson structures arise naturally in the realm of generalized geometry.
A holomorphic Poisson structure induces a deformation of the complex structure in a generalized sense,
whose cohomology
is obtained by twisting the Dolbeault $\dbar$-operator by the holomorphic Poisson bivector field.
 Therefore,
the cohomology space naturally appears as the limit of a spectral sequence of a double complex.
The first sheet of this spectral sequence is simply the Dolbeault cohomology with coefficients
in the exterior algebra of the holomorphic tangent bundle. We identify various necessary conditions
on compact complex manifolds on which this spectral sequence degenerates on the level of the second sheet.
The manifolds to our concern include all compact complex surfaces, K\"ahler manifolds,
and nilmanifolds with abelian complex structures or parallelizable complex structures.
\end{abstract}


\section{Introduction}

The algebraic geometry of Poisson brackets over complex manifolds were studied for some time \cite{Polish}.
In recent years, there are significant interests on holomorphic Poisson structures due to their emergence in
generalized complex geometry \cite{Marco} \cite{Hitchin-Generalized CY} \cite{Hitchin-Instanton}.
As such, one could consider
their deformations either as complex analytic objects
\cite{FM} \cite{Goto} \cite{Hitchin-holomorphic Poisson}, or as
generalized complex objects \cite{GPR}. On the other hand, there are classifications of holomorphic
Poisson structures on algebraic surfaces \cite{B-Marci}, and computations of their related cohomology
theory \cite{Hong-Xu}.

A holomorphic Poisson structure consists of a holomorphic bi-vector field $\Lambda$ on a complex manifold
such that $\lbra{\Lambda}{\Lambda}=0$. It is demonstrated in \cite{GPR} that the deformation of such structure as
a generalized complex structure is dictated by the cohomology
space $\oplus_kH^k_{\Lambda}$ of the differential operator
$\dbar_{\Lambda}=\dbar+\lbra{\Lambda}{-}$. This observation motivates the authors'
desire to find a way to compute the
cohomology spaces $H^k_{\Lambda}$. Since this space is naturally the limit of a spectral sequence of a double
complex due to the operators $\dbar$ and $\lbra{\Lambda}{-}$, it is an obvious question on when and how fast
this spectral sequence could degenerate. We will see that the first sheet of this spectral sequence is the
Dolbeault cohomology on the manifold with coefficients in the holomorphic polyvectors.
As this is a classical object, it becomes desirable to identify the conditions under
which the spectral sequence will degenerate at its second level.

In this paper, we identify three situations in which the spectral sequence of the
holomorphic Poisson double complex degenerates at its second level.
In the next chapter, we will explain the role of holomorphic Poisson structures in generalized geometry, and then
set up our computation of cohomology in the context of Lie bi-algebroids and their dual differentials
\cite{LWX}. Our reference for differential calculus of Lie algebroids is the book \cite{Mac}.
We will then set up the holomorphic Poisson double complex and its related spectral sequence.

We will make observation on three very different kinds of complex manifolds: complex parallelizable, nilmanifolds with abelian complex structures, and K\"ahlerian manifolds.

On complex parallelizable manifolds, we have the following result.

\

\noindent\textbf{Theorem 1.}\it
  Let $M$ be a complex parallelizable manifold with an invariant holomorphic Poisson structure $\Lambda$.
Let $\lieg^{1,0}$ be the vector space of holomorphic vector fields on $M$. Let
 $H^q_{\dbar}(\lieg^{*(0,1)})$ be the Chevalley-Eilenberg cohomology for the conjugation algebra
 $\lieg^{0,1}$, and $H_{\adL}^p(\lieg^{1,0})$ the cohomology of the operator $\adL=\lbra{\Lambda}{-}$
 on the exterior algebra
 of $\lieg^{1,0}$, then
 \[
 H^k_{\Lambda}(M)=\oplus_{p+q=k} H^q_{\dbar}(\lieg^{*(0,1)})\otimes H_{\adL}^p(\lieg^{1,0}).
 \]
 \rm

 This result will appear as Theorem \ref{parallelizable theorem} in this text. The proof of this result is not
 hard except for an application of Sakane's work \cite{Sakane} on Dolbeault cohomology on complex parallelizable
 manifolds because the holomorphic tangent bundle is trivial in this case. A typical example is the complex
 three-dimensional Iwasawa manifold.

The analysis of nilmanifolds takes significantly more effort.
Much are based on our experience with computation of various cohomology on nilmanifolds,
starting from Nomizu's work on de Rham cohomology \cite{Nomizu},
its generalization to Dolbeault cohomology \cite{Console-Fino} \cite{CFGU} and their
generalizations to cohomology with coefficients in holomorphic vector fields \cite{GPR} \cite{GMPP}
 \cite{MPPS-2-step}. In this paper, we need to deal with cohomology with coefficients in
holomorphic poly-vector fields. The result below generalizes what was known for Kodaira surface
\cite{Poon}.

\

\noindent\textbf{Theorem 2.} \it Suppose that $M=G/\Gamma$ is a 2-step
nilmanifold with an abelian complex structure.
Suppose that the center of $\lieg$ is real two-dimensional, then for
 any invariant holomorphic Poisson structure $\Lambda$,
 the Poisson spectral sequence degenerates on the second sheet.
\rm

\

This theorem will appear as Theorem \ref{theorem on nil}. In view of the work in \cite{CFP} and
\cite{Rolle}, we expect much room to generalize this theorem. There are also lots of examples.

Finally, to deal with K\"ahlerian manifolds, we adopt an idea from Hitchin in his paper
\cite{Hitchin-holomorphic Poisson} (see also \cite{Goto})  when he uses a contraction of the holomorphic
bivector field $\Lambda$ to pass a problem at hand to a consideration on Dolbeault cohomology.
This contraction yields a
homomorphism of double complex. Based on it, we prove the following result.

\

\noindent\textbf{Theorem 3.} \it
Suppose that $X$ is a holomorphic Poisson manifold with complex dimension $n$.
If the Fr\"{o}hlicher spectral sequence of its complex structure degenerates
at the $E_{1}$-level, then its holomorphic Poisson spectral sequence degenerates
at the $E_{n}$-level.
\rm

\

This theorem will appear as Theorem \ref{theorem through homo}.
Note that it is well known that the Fr\"{o}hlicher spectral sequence on a compact
complex surface always degenerates \cite{BVP}. The last theorem could be applied
to all compact complex surfaces as well as all compact K\"ahlerian manifolds \cite{Voisin}.

\section{Holomorphic Poisson Double Complex}

\subsection{Holomorphic Poisson Structures as Generalized Complex Structures}

A generalized complex manifold  \cite{Marco} \cite{Hitchin-Generalized CY} is a
smooth $(2n)$-dimensional manifold M equipped with a subbundle $L$ of the bundle
${\cal T}=(TM\oplus T^*M)_{\CC}$ such that
 \begin{itemize}
 \item  $L$ and its conjugate bundle $\overline{L}$ are transversal;
 \item  $L$ is maximally isotropic with respect to the natural pairing on $(TM\oplus T^*M)_{\CC}$;
  \item and the space of sections of $L$ is closed respective to the Courant bracket.
  \end{itemize}

 Equivalently, it is determined by a bundle automorphism $\cal J$ on $TM\oplus T^*M$
 such that ${\cal J}\circ {\cal J}=-$identity,  ${\cal J}$ is orthogonal with respect to the
  natural pairing, and the space of sections of
 the $+i$-eigenbundle with respect to ${\cal J}$ is closed with respect to the Courant bracket.

 Typical examples of generalized complex structures
 are classical complex structures and symplectic structures on a manifold.
 For a classical complex structure, the corresponding bundle $L$ as
 a generalized complex structure is $T^{1,0}M\oplus T^{*(0,1)}M$, the direct sum of the bundle of
 type (1,0)-vectors and the bundle of type (0,1)-forms.

Suppose that a manifold $M$ is equipped with a complex structure $J$.
A holomorphic Poisson structure on $(M, J)$ is a holomorphic bi-vector field $\Lambda$
such that with respect to the Schouten bracket, $\lbra{\Lambda}{\Lambda}=0$.
It follows that for any number $t$, the map
\begin{equation}
{\cal J}_{t\Lambda}=
\left(
\begin{array}
[c]{cc}
J & {\mbox {Im}}(t\Lambda) \\
0 & -J^*
\end{array}
\right)
\end{equation}
determines a generalized complex structure. Here by $J^*\alpha$ where $\alpha$ is a 1-form, we mean
$(J^*\alpha)X=\alpha(JX)$ for any vector field $X$.

At $t=0$, one recovers the generalized complex structure determined by the classical
complex structure $J$ on the manifold $M$. This is an example of deformation of generalized complex structures.

The corresponding bundles as generalized complex structure for ${\cal J}_{t\Lambda}$ is the pair of
bundles of graphes $L_{\overline{t\Lambda}}$ and ${\overline L}_{t\Lambda}$ where
\begin{equation}
{\overline L}_{t\Lambda}=\{{\overline\ell}+t\Lambda(\overline{\ell}): {\overline\ell}\in {\overline L}\}.
\end{equation}

\subsection{Holomorphic Poisson Cohomology}

Given a generalized complex structure,  the pair of bundles $L$
and $\overline L$ makes a Lie bi-algebroid. Via the canonical non-degenerate pairing on the bundle $\cal T$,
 the bundle $\overline L$ is complex linearly identified to the dual of $L$.
 Therefore, the Lie algebroid differential of $\overline L$ acts on $L$. It extends to a differential on
 the exterior algebra of $L$.
 For the calculus of Lie bialgebroids,
 we follow the conventions in \cite{Mac}. In particular, let $\rho$ be the anchor map on the Lie algebroid $L$.
 For any element $\Gamma$ in $C^\infty(M, \wedge^kL)$ and elements ${a}_1, \dots,
 {a}_{k+1}$ in $C^\infty(M, \wedge^k{\overline{L}})$, the differential  of $\Gamma$  is
 \begin{eqnarray*}
&& (\dbar_{L}\Gamma)({a}_1, \dots, {a}_{k+1})
=\sum_{r=1}^{k+1}(-1)^{r+1}\rho(a_r)(\Gamma({a}_1, \dots,{\hat{a}}_{r} ,\dots, {a}_{k+1}))\\
&\quad& \quad +
\sum_{r<s}(-1)^{r+s}\Gamma(\lbra{a_r}{a_s}, {a}_1, \dots,{\hat{a}}_{r} ,\dots,{\hat{a}}_{s} ,\dots,{a}_{k+1}).
\end{eqnarray*}

 For example, if $L$ comes from a complex structure,
 $L= T^{1,0}M\oplus T^{*(0,1)}M$. The differential is
 \begin{equation}
 \dbar: C^\infty(T^{1,0}M\oplus T^{*(0,1)}M)
 \to C^\infty((T^{1,0}M\otimes T^{*(0,1)}M ) \oplus T^{*(0,2)}M).
 \end{equation}
 Obviously, when $\oomega$ is a $(0,1)$-form, then $\dbar\oomega$ is the classical
 $\dbar$-operator in Dolbeault theory. When $X$ is a $(1,0)$-vector field,
 $\overline{Y}$ a $(0,1)$-vector field and $\omega$ a $(1,0)$-form, then
 \begin{equation}\label{dbar vector}
 (\dbar X)(\omega, \overline{Y}) = -{\overline{Y}}\omega(X)-X(\lbra{\omega}{\overline{Y}})
 = -{\overline{Y}}\omega(X)+({\cal{L}}_{\overline{Y}}\omega)(X)=\omega(\lbra{X}{\overline{Y}}).
 \end{equation}
 Therefore, $\dbar X$ is identical to the so-called Riemann-Cauchy operator as in \cite{Gau}.
 For future reference, we note the following observation, which could be proved
 by treating $\dbar$ as a Lie algebroid differential with respect to the
Lie algebroid $\overline{L}$ as we did in the calculus in equation (\ref{dbar vector}) above.

\begin{lemma}\label{dbar closed} For any $(2,0)$-field $\Lambda$, $(1,0)$-forms $\omega$, $\gamma$ and
$(0,1)$-vector field $\overline{Z}$,
\[
(\dbar\Lambda)(\omega, \gamma, \overline{Z})
=-\overline{Z}(\Lambda(\omega,\gamma))-\omega(\lbra{\overline{Z}}{\Lambda\gamma})
+\gamma(\lbra{\overline{Z}}{\Lambda\omega}).
\]
\end{lemma}

 \subsection{Poisson Double Complex}
We have seen that a holomorphic Poisson structure $\Lambda$
could be treated as a deformation of the complex structure
 $J$. While the pair of bundles
 $L_{\overline{\Lambda}}$ and ${\overline L}_{\Lambda}$ form a Lie bi-algebroid,
 so does the pair $L$ and ${\overline L}_{\Lambda}$ \cite{LWX}.

 \begin{proposition}\mbox{\rm\cite{GPR}} The Lie algebroid differential
 of ${\overline L}_{\Lambda}$ acting on $L$ is given by
 \begin{equation}
 \dbar_{\Lambda}=\dbar+\adL: C^\infty(T^{1,0}M\oplus T^{*(0,1)}M)\to
 C^\infty(\wedge^2(T^{1,0}M\oplus T^{*(0,1)}M)),
 \end{equation}
 where by definition $\adL(v+\overline{\omega})=\lbra{\Lambda}{v+\overline{\omega}}$ when
 $v$ is a smooth section of $T^{1,0}M$ and $\overline{\omega}$ is a smooth section of
 $T^{*(0,1)}M$.
 \end{proposition}

 The operator $\dbar_{\Lambda}$ extends to act on the exterior algebra of $T^{1,0}M\oplus T^{*(0,1)}M$:
 \begin{equation}
 \dbar_{\Lambda}: K^n\to K^{n+1}, \quad \mbox{ where } \quad
 K^n=C^\infty(\wedge^n(T^{1,0}M\oplus T^{*(0,1)}M)).
 \end{equation}
 Since $\Lambda$ is a holomorphic Poisson structure, $\dbar_{\Lambda}\circ \dbar_{\Lambda}=0$.
 Therefore, one has a complex with $\dbar_{\Lambda}$ being a differential.
 Starting with $k=0$, one has the $k$-th cohomology of this complex. It is denoted by
 $H^k_{\Lambda}(M)$. This cohomology space is called the
 $k$-th holomorphic Poisson cohomology of $\Lambda$; or simply Poisson cohomology.

 Since $\Lambda$ is holomorphic Poisson, it follows that
 \begin{equation}
 \dbar\circ \dbar=0, \quad \dbar \circ \adL+\adL\circ \dbar=0, \quad \adL\circ\adL=0.
 \end{equation}
Define $A^{p,q}=C^\infty(M, \wedge^pT^{1,0}M\otimes \wedge^qT^{*(0,1)}M)$, then
 \begin{equation}
 \dbar:A^{p,q}\to A^{p,q+1}; \quad \mbox{ and } \quad
 \adL: A^{p,q}\to A^{p+1,q}.
 \end{equation}
 Therefore, we obtain a double complex.
\[\begin{array}
[c]{ccccc}
&  &  &  & \\
& A^{p,q+1} & \longrightarrow & A^{p+1,q+1} & \\
\overline{\partial} & \uparrow &  & \uparrow & \\
& A^{p,q} & \longrightarrow & A^{p+1,q} & \\
&  & ad_{\Lambda} &  &
\end{array}
.
\]

For the double complex $\left(  A^{\ast,\ast},ad_{\Lambda},\overline{\partial}\right)$,
its associated single complex is $K^{n}=\oplus_{p+q=n}A^{p,q}$ and
the differential is
\[
\dbar_{\Lambda}=\overline{\partial}+\adL:K^{n}\longrightarrow K^{n+1}.
\]

 \begin{lemma} The triples $\left( A^{*,*}, \adL, \dbar\right)$ form a double complex. Its total cohomology is
 the holomorphic Poisson cohomology.
 \end{lemma}

The bi-grading yields two spectral sequences abutting to the Poisson cohomology.
We focus on the case when the filtration is given by
\begin{equation}
F^pK^n=\oplus_{p'+q=n, p'\geq p}A^{p',q}.
\end{equation}
Equivalently, the zero-th sheet of this spectral sequence $E_0^{p,q}$ is precisely $A^{p,q}$.
Since $\dbar_{\Lambda}=\dbar+\adL$ and $\adL$ does not raise the second index in the
bi-degree $(p,q)$, the restriction of the differential on the $0$-th sheet is
\begin{equation}
\dbar: A^{p,q}\to A^{p,q+1}.
\end{equation}
Therefore, the first sheet of this spectral sequence is the Dolbeault cohomology
\begin{equation}
 E_1^{p,q}=H^q(M, \Theta^p)
 \end{equation}
 where $\Theta^p$ is the $p$-th exterior power of the holomorphic tangent bundle of the complex manifold $M$.

 It follows that the next sheet of this spectral sequence is given by the cohomology of the differential
 \begin{equation}
 d_1^{p,q}=\adL: H^q(M, \Theta^p)\to H^q(M, \Theta^{p+1})
 \end{equation}
 for all $p+q\geq 0$, and the second sheet of this spectral sequence is
 \begin{equation}
 E_2^{p,q}=
 \frac{\mbox{kernel of } \adL: H^q(M, \Theta^p)\to
 H^q(M, \Theta^{p+1})}{\mbox{image of } \adL: H^q(M, \Theta^{p-1})
 \to H^q(M, \Theta^{p})}.
 \end{equation}

To facilitate further computation, we recall that the differential $d_2^{p,q}$ is computed as follows.
Suppose that \lq$v$' represents a class in $H^q(M, \Theta^p)$ and $\adL v$ represents a zero class in
$H^q(M, \Theta^{p+1})$, it represents an element in $E_2^{p,q}$.
 As $\adL v$ represents a zero class,
 there exists \lq$w$' a section of $T^{p+1,0}M\otimes T^{*(0,q-1)}M$ such that
$\dbar w=\adL v$. Since
\[
\dbar\adL w=-\adL \dbar w=-\adL\adL v=0,
\]
 $\adL w$ represents a class in
$E_2^{p+2, q-1}$. By definition, it is the image of the class of $v$ via the map $d_2^{p,q}$.

\subsection{Complex Projective Spaces as Examples}
Let $\Lambda$ be any non-trivial holomorphic section of the bundle $\Theta^2$ of the
complex projective space $M=\mathbb{CP}^n$. Assume that $\lbra{\Lambda}{\Lambda}=0$ so that
$\Lambda$ is a holomorphic Poisson structure. Since
$H^{q}(M,\Theta^{p})$ vanishes for all $q\geq1$, the first sheet of the
holomorphic Poisson double complex is reduced to its zero-th row:
\[
\begin{array}
[c]{ccccccc}
H^{0}(M, \mathcal{O}) & \overset{d_{1}^{0,0}=ad_{\Lambda}}{\longrightarrow} & H^{0}\left(
M, \Theta\right)  & \overset{d_{1}^{1,0}=ad_{\Lambda}}{\longrightarrow} & H^{0}\left(
M, \Theta^2 \right)  & \longrightarrow & \cdots
\end{array}.
\]
It is apparent that $d_{1}$ is not identically equal to zero. It is also
apparent that $d_{2}$ is identically equal to zero.

\subsection{Hopf Manifolds $M=S^1\times S^{2n-1}$}
Consider $\mathbb{C}^{n}$ with coordinates $z=(z_{1}, \dots, z_{n})$. Let $\lambda>1$
be a real number. It generates a one-parameter group of automorphism on
$\mathbb{C}^{n}$. The quotient of $\mathbb{C}^{n}\backslash\{0\}$ with respect
this group is diffeomorphic to the manifold  $M=S^1\times S^{2n-1}$. The complex
structure on $\mathbb{C}^{n}$ descends onto $M$ to define an integrable
complex structure. As a complex manifold, it is a principal elliptic fibration over the
complex projective space $N=\mathbb{CP}^{n-1}$:
\[
\pi: M=S^1\times S^{2n-1}\rightarrow N=\mathbb{CP}^{n-1}.
\]
Denote the fiber by $F$ and the holomorphic vertical vector field by $V$.
We have the following exact sequence of
holomorphic vector bundles over the Hopf manifold $M$:
\begin{equation}\label{basic exact}
0\to {\calo}_X(V)\to \Theta_M \to \pi^*\Theta_N\to 0,
\end{equation}
where ${\calo}_X(V)$ denote the vertical bundle of the projection $\pi$.
It is also a bundle trivialized by the vector
field $V$. Explicitly the vector field $V$ is simply the holomorphic vector
field generated by the group of dilations.
\begin{equation}
V=z_1\frac{\partial}{\partial z_1}+\cdots+z_n\frac{\partial}{\partial z_n}.
\end{equation}

For all $1\leq p\leq n-1$, the exact sequence (\ref{basic exact}) generates
\begin{equation}\label{exact}
0\to {\calo}_X(V)\otimes \pi^*\Theta_N^{p-1}\to \Theta_M^p \to \pi^*\Theta_N^{p}\to 0.
\end{equation}
We also have  ${\calo}_X(V)\otimes \pi^*\Theta_N^{n-1}\cong \Theta_M^{n}$.
Let $\oomega=\dbar \ln |z|^2$. It is a global $\dbar$-closed $(0,1)$-form on the Hopf manifold $M$.

\begin{proposition} On the Hopf manifold $M$,
\begin{itemize}
\item for all $p\geq 0$ and $q\geq 2$, $H^q(M, \Theta_M^p)=0$.
\item $H^0(M, \calo_M)=\langle 1\rangle$, and
$H^1(M, \calo_M)=\langle\oomega\rangle$.
\item $H^0(M, \Theta_M^n)=\langle V\rangle\otimes H^0(N, \Theta_N^{n-1})$,
$H^1(M, \Theta_M^n)=\langle \oomega\wedge V\rangle\otimes H^0(N, \Theta_N^{n-1})$.
\item
For all $1\leq p\leq n-1$,
\begin{eqnarray*}
H^0(M, \Theta_M^p) &=& H^0(N, \Theta_N^p)\oplus \langle V\rangle\otimes H^0(N, \Theta_N^{p-1}), \\
H^1(M, \Theta_M^p) &=& \langle \oomega\rangle\otimes H^0(M, \Theta_M^p)\\
&=&\langle \oomega\rangle\otimes H^0(N, \Theta_N^p)\oplus \langle \oomega\wedge V\rangle\otimes H^0(N, \Theta_N^{p-1}).
\end{eqnarray*}
\end{itemize}
\end{proposition}
The proof of this proposition is a straight-forward application of the
Leray spectral sequence with respect to the
fiberation $\pi$. Similar computation could be found in \cite{GMPP} and \cite{Hoffer}.
Below we outline a proof. First of all, by virtue of the dimension of the fiber and the global $(0,1)$-form $\oomega$
on the manifold $M$, one proves that the direct image sheaves of $R^q\pi_*\calo_M={0}$ for $q\geq 2$, and
\[
R^0\pi_*\calo_M\cong \calo_N, \quad R^1\pi_*\calo_M\cong \calo_N(\oomega).
\]
By the Projective Formula, for all $k\geq 0$, $R^q\pi_*\pi^*\Theta^k_N={0}$ for $q\geq 2$, and
\[
R^0\pi_*\pi^*\Theta^k_N\cong \Theta^k_N, \quad R^1\pi_*\pi^*\Theta^k_N\cong \calo_N(\oomega)\otimes \Theta^k_N.
\]
As the second sheet of the Leray spectral sequence for $\pi^*\Theta^k_N$ is
\[
E_2^{p,q}(\pi^*\Theta^k_N)=H^p(N, R^q\pi_*\pi^*\Theta^k_N).
\]
$E_2^{p,q}=0$ for all $q\geq 2$. Moreover,
\[
E_2^{p,0}(\pi^*\Theta^k_N)=H^p(N, \Theta^k_N), \quad
E_2^{p,1}(\pi^*\Theta^k_N)=\langle\oomega\rangle\otimes H^p(N, \Theta^k_N).
\]
Recall that $N$ is the complex projective space $\mathbb{CP}^{n-1}$. Therefore, the above cohomology
spaces vanish except when $p=0$. It follows immediately the differential $d_2^{p,q}$ is identically zero,
and hence the Leray
spectral sequence degenerates to $\oplus_{p+q=m}E^{p,q}_2(\pi^*\Theta^k_N)=H^{m}(M, \pi^*\Theta^k_N)$.
In particular,
\begin{equation}
H^{0}(M, \pi^*\Theta^k_N)\cong H^0(N, \Theta^k_N), \quad
H^{1}(M, \pi^*\Theta^k_N)\cong \langle\oomega\rangle\otimes H^0(N, \Theta^k_N).
\end{equation}
By the exact sequence (\ref{exact}), one concludes immediately, that for all $p\geq 0$ and $q\geq 2$,
$H^q(M, \Theta_M^p)=0$. By chasing diagram, one could actually show that the induced long exact sequence of
(\ref{exact}) splits. The claim in the proposition above follows.

\

Given the proposition above, one notes that the non-zero terms of the first
sheet of the Poisson double complex with respect to any holomorphic
Poisson structure on the Hopf manifold are contained in the zero-th and the
first rows only.
To compute the differential on the second sheet explicitly, note that the first sheet is given by
\[
\begin{array}{cccccccc}
\left\langle \overline{\omega }\right\rangle & \overset{ad_{\Lambda }}{
\rightarrow } & \left\langle \overline{\omega }\right\rangle \otimes
H^{0}(M,\Theta_M ) & \overset{ad_{\Lambda }}{\rightarrow } & \left\langle
\overline{\omega_M }\right\rangle \otimes H^{0}(M,\Theta_M^{2}) & \overset{
ad_{\Lambda }}{\rightarrow } & \left\langle \overline{\omega }\right\rangle
\otimes H^{0}(M,\Theta_M^{3}) & \overset{ad_{\Lambda }}{\rightarrow } \\
&  &  &  &  &  &  &  \\
\left\langle 1\right\rangle & \overset{ad_{\Lambda }}{\rightarrow } &
H^{0}(M,\Theta_M ) & \overset{ad_{\Lambda }}{\rightarrow } & H^{0}(M,\Theta_M
^{2}) & \overset{ad_{\Lambda }}{\rightarrow } & H^{0}(M,\Theta_M^{3}) &
\overset{ad_{\Lambda }}{\rightarrow }
\end{array}.
\]

To prepare our computation, note that when we treat $\mathbb{CP}^{n-1}$ as a
homogeneous space, then its space of holomorphic vector fields is
$\mathfrak{sl}(n,\mathbb{C})$. They are generated by linear action of $SU(n)$ on
$\mathbb{C}^{n}.$

Let $A\in \mathfrak{sl}(n,\mathbb{C})$, and define on $\mathbb{C}^n$  the function
\[
f_A(z)=\frac{1}{\left\vert z\right\vert ^{2}}\overline{z}
^{T}Az.
\]
By homogeneity, this function is well-defined on the Hopf manifold $M$.
It is straightforward  to show that
 $\mathcal{L}_{A}
\overline{\omega }=\overline{\partial }f_{A}$ and
$\mathcal{L}_{V}f_{A}=0.$

\

Consider the bi-vector field $\Lambda =V\wedge A.$ Since $A$ and $V$ commute,
$\Lambda$ is a holomorphic Poisson structure.
We claim that its associated Poisson double complex
degenerates at the $E_{2}$-level.

Since the only non-trivial elements in the first sheet of the double complex are in the zero-th
and the first row, $E_{2}^{p,q}$ vanishes except possibly when $
(p,q)=(p,0)$ and $(p,q)=(p,1).$ Therefore, the only possibly non-trivial
differential on the second sheet are $d_{2}^{p,1}$. According to the
structure of the cohomology spaces, each term in $E_{2}^{p,1}$ is
represented by an element of the form $\overline{\omega }\wedge \Gamma $ for
some $\Gamma $ in $H^{0}(M,\Theta ^{p})$. Since $\call_V(\oomega)=0$ and $\call_V(\Gamma)=0$,
\begin{eqnarray*}
  &&
\lbra{V\wedge A}{\overline{\omega }\wedge \Gamma}
=V\wedge \lbra{A}{\overline{\omega }\wedge \Gamma}
=V\wedge \mathcal{L}_{A}\overline{\omega }\wedge \Gamma
+(-1)^{p}V\wedge\lbra{A}{\Gamma}\wedge \overline{\omega } \\
&=&V\wedge \overline{\partial }f_{A}\wedge \Gamma
+(-1)^{p}V\wedge \lbra{A}{\Gamma} \wedge \overline{\omega }
=-\overline{\partial }\left( f_{A}V\wedge \Gamma \right)
+(-1)^{p}V\wedge \lbra{A}{\Gamma} \wedge \overline{\omega }.
\end{eqnarray*}
So, ${\adL}({\overline{\omega }\wedge \Gamma})$
represents the zero class in
$\left\langle \overline{\omega }\right\rangle\otimes H^{0}(M,\Theta ^{p+1})$
if and only if
$\lbra{A}{\Gamma}$ is equal to zero.
Under this condition,
${\adL}({\overline{\omega }\wedge \Gamma} )
=-\overline{\partial }\left( f_{A}V\wedge \Gamma \right)$, and
 $d_{2}^{p,1}\left( \overline{\omega }\wedge \Gamma \right) $ is
represented by
\begin{eqnarray*}
&&-{\adL}({f_{A}V\wedge \Gamma})
=-\lbra{ V\wedge A}{f_{A}V\wedge \Gamma}=A\wedge \lbra{V}{f_{A}V\wedge \Gamma}
-V\wedge \lbra{A}{f_{A}V\wedge \Gamma }  \\
&=&A\wedge (\mathcal{L}_{V}f_{A}) V\wedge \Gamma
+A\wedge f_{A}\lbra{V}{V\wedge \Gamma} -V\wedge (\mathcal{L}_{A}f_{A}) V\wedge \Gamma
-V\wedge f_{A}\lbra{A}{V\wedge \Gamma}.
\end{eqnarray*}
Note that the first summand is equal to zero because $\mathcal{L}_{V}f_{A}=0$.
The second summand is equal to zero because $V$ commutes with
every element in $H^{0}(M,\Theta ^{p})$. The third summand is equal to zero
by skew-symmetry of exterior multiplication. The last summand is equal to
zero because as noted above when the term
$\lbra{\adL}{\overline{\omega }\wedge \Gamma }$ represents the zero class in
$H^{1}(M,\Theta^{p+1})$,  $\lbra{A}{\Gamma}=0.$ Therefore, $d_{2}^{p,1}= 0$ for all $p$.

\bigskip

In the rest of this article, we investigate the degeneracy of
holomorphic Poisson spectral sequence in various special settings.

\section{Nilmanifolds}

A compact manifold $M$ is called a nilmanifold if there exists a simply-connected nilpotent Lie group $G$
and a lattice subgroup $\Gamma$ such that $M$ is diffeomorphic to $G/\Gamma$. We denote the Lie algebra of
$G$ by $\lieg$, and the center of $\lieg$ by $Z(\lieg)$. The step of the
nilmanifold is the nilpotence of the Lie algebra $\lieg$. A left-invariant complex
structure $J$ on $G$ is said to be abelian if on the Lie algebra $\lieg$, it satisfies the conditions
$J\circ J=-$identity and $\lbra{JA}{JB}=\lbra{A}{B}$ for all $A$ and $B$ in the Lie algebra $\lieg$.
If one complexifies the algebra $\lieg$ and denotes the $+i$ and $-i$ eigen-spaces of $J$ respectively
by $\lieg^{1,0}$ and $\lieg^{0,1}$, then the invariant complex structure $J$ being abelian is equivalent
to the complex algebra $\lieg^{1,0}$ being abelian.

Denote $\wedge^k\lieg^{1,0}$ and $\wedge^k{\lieg}^{*(0,1)}$ respectively by
$\lieg^{k,0}$ and ${\lieg}^{*(0,k)}$.

\subsection{Cohomology of Nilmanifolds}

Now assume that the Lie algebra $\lieg$ is 2-step nilpotent, i.e. $[\lieg,\lieg]\subset Z(\lieg)$.
In such a case, we call $M=G/\Gamma$ a 2-step nilmanifold \cite{MPPS-2-step}.

On the nilmanifold $M$, we consider $\lieg^{k,0}$ as invariant $(k,0)$-vector fields, and
${\lieg}^{*(0,k)}$ as invariant $(0,k)$-forms. It yields an inclusion map
\[
{\lie g}^{p,0}\otimes {\lie g}^{*(0,q)}
\hookrightarrow C^\infty(M, \wedge^pT^{1,0}M\otimes \wedge^qT^{*(0,1)}M).
\]
When the complex structure is also invariant, $\dbar$ sends
${\lie g}^{p,0}\otimes {\lie g}^{*(0,q)}$ to
${\lie g}^{p,0}\otimes {\lie g}^{*(0,q+1)}$. The concerned cohomology is
\begin{equation}\label{hqp}
H^q( \lieg^{p,0})=
\frac
{\mbox{kernel of } \dbar: {\lieg}^{*(0,q)}\otimes {\lieg}^{p,0}\to {\lieg}^{*(0,q+1)}\otimes {\lieg}^{p,0}}
{\mbox{image of } \dbar: {\lieg}^{*(0,q-1)}\otimes {\lieg}^{p,0}\to {\lieg}^{*(0,q)}\otimes {\lieg}^{p,0}}
.\end{equation}
The inclusion map yields an inclusion of cohomology:
\[
H^q({\lieg}^{p,0})\hookrightarrow H^q(M, \Theta^p).
\]

\begin{theorem}\label{quasi isomorphic} On a 2-step nilmanifold $M$ with an invariant
abelian complex structure, the inclusion ${\lie g}^{p,0}\otimes {\lie g}^{*(0,q)}$ in
$C^\infty(M, \wedge^pT^{1,0}M\otimes \wedge^qT^{*(0,1)}M)$ induces an isomorphism of cohomology. In other words,
\[
H^q({\lieg}^{p,0})\cong H^q(M, \Theta^p).
\]
\end{theorem}

We need some preparations before we give a proof for this theorem.

Let $C$ be the center of $G$ and $\psi: G\rightarrow G/C$ the quotient map. Since $G$ is 2-step nilpotent,
$G/C$ is abelian. Consider $M=G/\Gamma$ and $N=\psi(G)/\psi(\Gamma)$.
We have a holomorphic fibration $\Psi: M\rightarrow N$ whose fiber is isomorphic to $F=C/(C\cap \Gamma)$.

 Let $\liec=Z(\lieg)$, $\liet=\lieg/Z(\lieg)$.
 Below are some facts shown in Sections 2 and 3 of \cite{MPPS-2-step}.
 Since $\lieg$ is 2-step nilpotent, $\liet$ is abelian.
 As a vector space, $\lieg^{1,0}=\liet^{1,0}\oplus \liec^{1,0}$,
 and  $\lieg^{*(1,0)}=\liet^{*(1,0)}\oplus \liec^{*(1,0)}$.
 The only nontrivial Lie brackets in $\lieg^{1,0}\oplus \lieg^{0,1}$
 are of the form
 $ [\liet^{1,0},\liet^{0,1}]\subset \liec^{1,0}\oplus \liec^{0,1}.$

Explicitly, there exists a real basis $\{X_k,JX_k:1\leq k\leq n\}$ for $\liet$ and
$\{Z_\ell, JZ_\ell:1\leq \ell \leq m\}$ a real basis for $\liec$. The corresponding complex bases
for $\liet^{1,0}$ and $\liec^{1,0}$ are respectively composed of the following elements:
\begin{equation}
T_k=\frac{1}{2}(X_k-iJX_k), \quad \mbox{ and } \quad W_\ell=\frac{1}{2}(Z_\ell-iJZ_\ell).
\end{equation}
The structure equations of $\lieg$ are determined by
\begin{equation}\label{structure eq}
\lbra{\overline{T}_k}{T_j}=\sum_{\ell}E_{kj}^\ell W_\ell-\sum_\ell{\overline{E}}_{jk}^\ell{\overline{W}}_\ell
\end{equation}
for some constants $E_{kj}^\ell$.
Let $\{\omega^k: 1\leq k\leq n\}$ be the dual basis for $\liet^{*(1,0)}$, and let
$\{\rho^\ell: 1\leq \ell\leq m\}$ be the dual basis for $\liec^{*(1,0)}$.
By formula (\ref{dbar vector}),
\begin{equation}
\dbar {T_j}=\sum_{k, \ell}E_{kj}^\ell\oomega^k\wedge W_\ell,
\end{equation}
The dual structure equations for (\ref{structure eq}) are
\begin{equation}\label{dual 1}
d\rho^\ell=\sum_{i,j}E_{ji}^\ell\omega^i\wedge\oomega^j.
\end{equation}
Equivalently,
\begin{equation}
d\overline{\rho}^\ell=-\sum_{i,j}\overline{E}^\ell_{ji}\omega^j\wedge\oomega^i.
\end{equation}
It follows that
\begin{equation}\label{adj-T-on-rho-bar}
\lbra{T_j}{\overline{\rho}^\ell}={\cal L}_{T_j}{\overline{\rho}^\ell}=-\sum_i\overline{E}^\ell_{ji}\oomega^i.
\end{equation}
In summary, the  operator $d=\partial+\bar{\partial}$ satisfies
 \begin{equation}\label{Eqt:dliegstar}
 \partial\liet^{*(1,0)}=0,\qquad \bar{\partial}\liet^{*(1,0)}=0,\qquad
\partial\liec^{*(1,0)}=0,
\qquad \bar{\partial}\liec^{*(1,0)}\subset \liet^{*(1,1)}.
   \end{equation}
Moreover, the co-boundary operator
$\bar{\partial}: \lieg^{1,0}\rightarrow \lieg^{*(0,1)}\otimes \lieg^{1,0}$ satisfies
 \begin{equation}\label{Eqt:barpartiallieg10}
  \bar{\partial}\liec^{1,0}=0,\qquad
\bar{\partial}\liet^{1,0} \subset \liet^{*(0,1)}\otimes \liec^{1,0}.
   \end{equation}

\

\begin{lemma} Let $\Theta_N$ be the tangent sheaf of $N$, the quotient manifold. Then
$R^q\Psi_*(\wedge^l \Psi^*\Theta_N)=\liec^{*(0,q)}\otimes \wedge^l\Theta_N.$
\end{lemma}
 \bproof By Lemma 4 in \cite{MPPS-2-step},
 we have $R^q\Psi_*(\mathcal{O}_M)= \liec^{*(0,q)}\otimes \mathcal{O}_N$.
 Thus, by the projective formula,
\[
R^q\Psi_*(\wedge^l \Psi^*\Theta_N)=R^q\Psi_*( \Psi^*\wedge^l\Theta_N)
 =R^q\Psi_*(\mathcal{O}_M)\otimes \wedge^l\Theta_N=\liec^{*(0,q)}\otimes \wedge^l\Theta_N.
 \]
 \eproof

\begin{lemma}\label{Lemma:intermiedietfacts}
$H^k(N,\wedge^l \Theta_N)=\liet^{*(0,k)}\otimes \liet^{l,0},$ and
$H^k(M,\wedge^l\Psi^* \Theta_N )=\lieg^{*(0,k)}\otimes \liet^{l,0}$.
\end{lemma}
\bproof The proof is essentially the same as that of Lemma 5 in \cite{MPPS-2-step}.
Here we give a sketch. First note that $N$ is an abelian variety. It follows that
 $\Theta_N$ is holomorphically trivial, and is isomorphic to $\liet^{1,0}$.
Therefore, $H^k(N,\wedge^l \Theta_N)=\liet^{*(0,k)}\otimes \liet^{l,0}$.

To compute $H^k(M,\wedge^l\Psi^* \Theta_N )$, we use the standard
Leray spectral sequence for the fibration $\Psi$. By the previous lemma, the second sheet is given by
\begin{align*}
E^{p,q}_2&=H^p(N,R^q\Psi_*(\wedge^l\Psi^*\Theta_N))=H^p(N,\liec^{*(0,q)}\otimes \wedge^l\Theta_N)\\
&=\liec^{*(0,q)}\otimes H^p(N, \wedge^l\Theta_N)=\liec^{*(0,q)}\otimes \liet^{*(0,p)}\otimes \liet^{l,0}.
\end{align*}
Thus $d_2: E^{p,q}_2\rightarrow E^{p+2,q-1}_2$ is a map
$$
\liec^{*(0,q)}\otimes \liet^{*(0,p)}\otimes \liet^{l,0}\rightarrow
\liec^{*(0,q-1)}\otimes \liet^{*(0,p+2)}\otimes \liet^{l,0}.
$$
But $d_2$ is essentially given by $\bar{\partial}$.
Equalities \eqref{Eqt:dliegstar} and \eqref{Eqt:barpartiallieg10} imply that
$d_2=0$. It follows that
\begin{align*}
H^k(M,\wedge^l\Psi^* \Theta_N )&=\oplus_{p+q=k}  E^{p,q}_2
=\oplus_{p+q=k}  \liec^{*(0,q)}\otimes \liet^{*(0,p)}\otimes \liet^{l,0}\\
&=\lieg^{*(0,k)}\otimes \liet^{l,0}.
\end{align*}
\eproof
\begin{lemma}\label{Lemma:filteration}
Suppose that $$0\rightarrow Z\rightarrow E\stackrel{\rho}{\longrightarrow} Q\rightarrow 0$$
is an exact sequence of holomorphic vector bundles. Then
there induces a series of exact sequences
\begin{eqnarray*}
&0\rightarrow E^{(1)}\rightarrow \wedge^l E \rightarrow \wedge^l Q;&\\
&0\rightarrow E^{(2)}\rightarrow E^{(1)} \rightarrow  \wedge^{l-1} Q\otimes Z;&\\
&\cdots\cdots\cdots&\\
&0\rightarrow E^{(r+1)}\rightarrow E^{(r)} \rightarrow    \wedge^{l-r} Q \otimes \wedge^r Z;&\\
&\cdots\cdots\cdots&\\
&0\rightarrow E^{(l)}\rightarrow E^{(l-1)} \rightarrow Q\otimes \wedge^{l-1} Z ,&
\end{eqnarray*}
where $E^{(l)}=\wedge^l Z$.
\end{lemma}
\bproof For any $m\geq 0$ and $n\geq 0$, we define a derivative
 operator $D_{\rho}:~\wedge^m Q\otimes \wedge^n E\rightarrow \wedge^{m+1} Q\otimes \wedge^{n-1} E$ by
\begin{align*}
&D_{\rho}(  q_1\wedge\cdots\wedge q_m \otimes e_1 \wedge\cdots\wedge e_n)\\
=&\sum_{i=1}^n (-1)^{i-1} q_1\wedge\cdots\wedge q_m
\wedge \rho(e_i)\otimes  e_1\wedge\cdots {\widehat{e_i}}\cdots\wedge e_n .
\end{align*}
Now consider $$D_{\rho}^l: \wedge^l E\rightarrow \wedge^l Q.$$
Let $E^{(1)}=\mathrm{Ker}D_{\rho}^l$.  Similarly, for each $r=1,\cdots, l$, we consider
$$D_{\rho}^r: \wedge^{l} E\rightarrow  \wedge^{l+r} Q \otimes \wedge^{l-r} E, $$
and let $E^{(l-r+1)}=\mathrm{Ker}D_{\rho}^r$.

It is obvious that $\wedge^l Z=E^{(l)}\subset E^{(l-1)}\subset \cdots \subset E^{(1)}$.
The fact that $E^{(r)}/E^{(r+1)}\cong   \wedge^{l-r} Q \otimes \wedge^r Z$ is also easily seen.
\eproof

\newcommand{\huaF}{\mathcal{F}}
\newcommand{\partialbar}{\bar{\partial}}

\bproof[of Theorem \ref{quasi isomorphic}] Let $E=\Theta_M$,
$Z=\mathcal{O}_M\otimes \liec^{1,0}  $, $Q=\Psi^*\Theta_N$.
Then we have the following exact sequence of holomorphic vector bundles:
$$0\rightarrow Z\rightarrow E\stackrel{\rho}{\longrightarrow} Q\rightarrow 0.$$
By Lemma \ref{Lemma:filteration}, we have a filtration of $\wedge^l E=\wedge^l \Theta_M$:
\[
\mathcal{O}_M\otimes \liec^{l,0}=\wedge^l Z =E^{(l)}\subset E^{(l-1)}
\subset \cdots \subset E^{(1)}\subset E^{(0)}=\wedge^l E=\wedge^l \Theta_M.
\]

Moreover, the associated graded spaces are:
\begin{eqnarray*}
&G^0= E^{(0)} / E^{(1)}\cong \wedge^l Q=\Psi^*(\wedge^{l} \Theta_N);& \\
&G^1=E^{(1)} / E^{(2)}\cong  \wedge^{l-1} Q\otimes Z=\Psi^*(\wedge^{l-1} \Theta_N)\otimes \liec^{1,0};&\\
&\cdots\cdots\cdots&\\
&G^r= E^{(r)} / E^{(r+1)}=\wedge^{l-r} Q\otimes \wedge^r Z= \Psi^*(\wedge^{l-r} \Theta_N)\otimes \liec^{r,0};&\\
&\cdots\cdots\cdots&\\
&G^{l-1}=E^{(l-1)} / E^{(l)}\cong  Q\otimes \wedge^{l-1}Z = \Psi^*(  \Theta_N)\otimes \liec^{l-1,0}.&
\end{eqnarray*}

Accordingly, we have a filtration of the co-chain complex
$C^\bullet = T^{*(0,\bullet)}\otimes \wedge^l \Theta_M$:
$$  T^{*(0,\bullet)}\otimes \liec^{l,0}
=\huaF^{(l)} C^{\bullet}\subset \huaF^{(l-1)} C^{\bullet} \subset \cdots
\subset\huaF^{(1)} C^{\bullet}\subset \huaF^{(0)} C^{\bullet}= C^\bullet,
$$
where $\huaF^{(r)}C^\bullet=T^{*(0,\bullet)}\otimes E^{(r)}$.
Thus there associates a spectral sequence: $E_r^{p,q}$, which starts with
$$
E^{p,q}_0= \huaF^{(p)}C^{p+q}/ \huaF^{(p+1)}C^{p+q}=
T^{*(0,\bullet)}\otimes G^p\cong T^{*(0,\bullet)}\otimes \Psi^*(\wedge^{l-p} \Theta_N)\otimes \liec^{p,0}
$$
and $d_0=\bar{\partial}$. It follows from Lemma
\ref{Lemma:intermiedietfacts} that we have $$ E^{p,q}_1=
H^{p+q}(G^p)=H^{p+q}( \Psi^*(\wedge^{l-p} \Theta_N)\otimes
\liec^{p,0})=\lieg^{*(0,p+q)}\otimes\liet^{l-p,0}\otimes\liec^{p,0}.
$$ The associated $d_1: E^{p,q}_1\rightarrow E^{p+1,q}_1$ should be exactly $$\bar{\partial}:
\lieg^{*(0,p+q)}\otimes\liet^{(l-p,0)}\otimes\liec^{(p,0)}
\rightarrow
\lieg^{*(0,p+1+q)}\otimes\liet^{(l-p-1,0)}\otimes\liec^{(p+1,0)} .$$
Here we use the fact that $\lieg$ is 2-step nilpotent, namely the
relations in \eqref{Eqt:dliegstar} and \eqref{Eqt:barpartiallieg10}.
Hence
$$E^{p,q}_2= H^{p+q}_{d_1}(E^{p,q}\rightarrow E^{p+1,q})=
\frac{\text{Ker}(\bar{\partial})\cap
\lieg^{*(0,p+q)}\otimes\liet^{(l-p,0)}\otimes\liec^{(p,0)}
}{\bar{\partial}(\lieg^{*(0,p-1+q)}\otimes\liet^{(l-p+1,0)}\otimes\liec^{(p-1,0)}
)}.
$$

The operator $d_2: E^{p,q}_2\rightarrow E^{p+2,q-1}_2$ is
essentially that of $\bar{\partial}$. So the previous expression
implies that $d_2=0$. So  $E^{p,q}_2=E^{p,q}_{\infty}$. Therefore,
one has
\begin{eqnarray*}
&& H^k(\wedge^l \Theta_M)=H^k_{\partialbar}(C^\bullet)\cong \bigoplus_{p+q=k} E^{p,q}_2\\
&=&\bigoplus_{p+q=k} \frac{\text{Ker}(\partialbar)\cap
\lieg^{*(0,p+q)}\otimes\liet^{l-p,0}\otimes\liec^{p,0} }{\partialbar(\lieg^{*(0,p-1+q)}
\otimes\liet^{l-p+1,0}\otimes\liec^{p-1,0}  )}
=\bigoplus_{p} \frac{\text{Ker}(\partialbar)
\cap \lieg^{*(0,k)}\otimes\liet^{l-p,0}
\otimes\liec^{p,0} }{\partialbar(\lieg^{*(0,k-1 )}\otimes\liet^{l-p+1,0}\otimes\liec^{p-1,0}  )}\\
& =&  \frac{\text{Ker}(\partialbar)\cap
(\bigoplus_{p}\lieg^{*(0,k)}\otimes\liet^{l-p,0}\otimes\liec^{p,0}
)}{\bigoplus_{p} \partialbar(\lieg^{*(0,k-1
)}\otimes\liet^{l-p+1,0}\otimes\liec^{p-1,0}  )}
= \frac{\text{Ker}(\partialbar)\cap \lieg^{*(0,k)}\otimes
\lieg^{l,0} }{\partialbar(\lieg^{*(0,k-1 )}\otimes \lieg^{l,0}  )}\\
&=& H^k_{\partialbar}(\lieg^{l,0}),
\end{eqnarray*}
as required.
\eproof

\subsection{Degeneracy of Holomorphic Poisson Double Complex}
As a consequence of Theorem \ref{quasi isomorphic},
\begin{equation}
 E_2^{p,q}=
 \frac{\mbox{kernel of } \adL:  H^q({\lieg}^{p,0})
 \to  H^q({\lieg}^{p+1,0})}{\mbox{image of }
 \adL:  H^q({\lieg}^{p-1,0}) \to
  H^q({\lieg}^{p,0})}
 \end{equation}
where $H^q( \lieg^{p,0})$ is given in (\ref{hqp}).

As an application, we consider the case when the nilmanifold is 2-step with real two-dimensional center.
Since
\[
\lieg^{1,0}=\liet^{1,0}\oplus \liec^{1,0} \quad \mbox{ and }
\quad \lieg^{*(0,1)}=\liet^{*(0,1)}\oplus \liec^{*(0,1)},
\]
and the complex dimension of $\liec^{1,0}$ is equal to one,
\[
\lieg^{p,0}=\liet^{p,0}\oplus \liet^{p-1,0}\otimes \liec^{1,0},
\quad \quad
\lieg^{*(0,q)}=\liet^{*(0,q)}\oplus \liet^{*(0,q-1)}\otimes\liec^{*(0,1)}.
\]

\begin{lemma} \label{image of dbar}
$
\lieg^{*(0,1)}\oplus \liec^{1,0}\subseteq \ker\dbar$ and
$\dbar \liet^{1,0}\subseteq \liet^{*(0,1)}\otimes \liec^{1,0}$.
Moreover,
\begin{equation}
\dbar (\liet^{*(0,a)}\otimes \liec^{*(0,b)}\otimes \liet^{k,0}\otimes \liec^{\ell,0})
\subseteq \liet^{*(0,a+1)}\otimes \liec^{*(0,b)}\otimes \liet^{k-1,0}\otimes \liec^{\ell+1,0}.
\end{equation}
\end{lemma}

In subsequent presentation, we suppress the notations for the vector spaces,
and simply keep track of the  quadruple of indices to indicate the components.
With this notation, the statement above is summarized as
\begin{equation}\label{dbar action}
\dbar(a, b; k, \ell)\subseteq (a+1, b; k-1, \ell+1).
\end{equation}
Due to dimensional assumption on $\liec^{1,0}$, when $\ell\geq 1$, then the space
$(a+1, b; k-1, \ell+1)$ is trivial.

Next, we turn our attention to the effect of an invariant holomorphic Poisson structure.
Given the dimension restriction on $\liec^{1,0}$, $\lieg^{2,0}=\liec^{1,0}\otimes \liet^{1,0}\oplus
\liet^{2,0}$. With regard to this decomposition, any $\Lambda$ in $\lieg^{2,0}$ is the sum of
$\Lambda_1$ and $\Lambda_2$ where,
$\Lambda_1=W\wedge T$ with $W\in \liec^{1,0}$,  $T\in\liet^{1,0}$, and
$\Lambda_2$ is in $\liet^{2,0}$.

Since the complex structure is abelian,
$\lbra{\Lambda}{\lieg^{1,0}}=0$. One could also see
from the structure equations (\ref{dual 1}) that $\lbra{\Lambda}{\oomega^k}=0$.
By (\ref{adj-T-on-rho-bar}),
\begin{equation}
\lbra{\Lambda_1}{\overline{\rho}}={W}\wedge\lbra{T}{\overline{\rho}}\subseteq
\liet^{*(0,1)}\otimes \liec^{1,0}.
\end{equation}
\begin{lemma}\label{image of Lambda}
 When $\Lambda=\Lambda_1+\Lambda_2$, then
 $\adL(\lieg^{1,0})=0,$ $\adL(\liet^{*(0,1)})=0,$ and
\begin{equation}
{\mbox\rm ad}_{\Lambda_1}(\liec^{*(0,1)})\subseteq
\liet^{*(0,1)}\otimes \liec^{1,0}, \quad
{\mbox\rm ad}_{\Lambda_2}(\liec^{*(0,1)})\subseteq \liet^{*(0,1)}\otimes \liet^{1,0}.
\end{equation}
\end{lemma}

As a consequence of the above lemma,
$
{\mbox\rm ad}_{\Lambda_1} (\liet^{*(0,a)}\otimes \liec^{*(0,b)}\otimes \liet^{k,0}\otimes \liec^{\ell,0})$ is
contained in $\liet^{*(0,a+1)}\otimes \liec^{*(0,b-1)}\otimes \liet^{k,0}\otimes \liec^{\ell+1,0},
$
and
$
{\mbox\rm ad}_{\Lambda_2} (\liet^{*(0,a)}\otimes \liec^{*(0,b)}\otimes \liet^{k,0}\otimes \liec^{\ell,0})$ is
contained in $\liet^{*(0,a+1)}\otimes \liec^{*(0,b-1)}\otimes \liet^{k+1,0}\otimes \liec^{\ell,0}.$
Using our shorthand notations, we have
\begin{equation}\label{adL action}
{\mbox\rm ad}_{\Lambda_1}(a,b;k,\ell)=(a+1, b-1; k, \ell+1), \quad
{\mbox\rm ad}_{\Lambda_2}(a,b;k,\ell)=(a+1, b-1; k+1, \ell).
\end{equation}
As the range of ${\mbox\rm ad}_{\Lambda_1}$ and ${\mbox\rm ad}_{\Lambda_2}$ are in different spaces,
\[
\ker \adL=\ker {\mbox\rm ad}_{\Lambda_1}\cap \ker {\mbox\rm ad}_{\Lambda_2}.
\]

\begin{theorem}\label{theorem on nil} Suppose that $M=G/\Gamma$ is a 2-step
nilmanifold with an abelian complex structure.
Suppose that the center of $\lieg$ is real two-dimensional, then for
 any invariant holomorphic Poisson structure $\Lambda$,
 the Poisson spectral sequence degenerates on the second sheet.
\end{theorem}

We need to prove that $d_2^{p,q}\equiv 0$ for all $p,q$. Given the dimension assumption on $\liec^{1,0}$,
Using our shorthand notations, the space $\lieg^{*(0,q)}\otimes \lieg^{p,0}$ decomposes into four components:
\[
(q,0;p,0); \quad (q-1,1;p,0); \quad (q,0;p-1,1); \quad (q-1,1;p-1,1).
\]
Suppose that $A$ is in $\lieg^{*(0,q)}\otimes \lieg^{p,0}$ and $\dbar A=0$.
Consider $\adL(A)={\mbox\rm ad}_{\Lambda_1}A+{\mbox\rm ad}_{\Lambda_2}A$.
By (\ref{adL action}),
$(q,0;p,0)$, $(q,0;p-1,1)$ and  $(q-1,1;p-1,1)$ are in $\ker {\mbox\rm ad}_{\Lambda_1}$, and
\[
{\mbox\rm ad}_{\Lambda_1}(A)\in {\mbox\rm ad}_{\Lambda_1}\lieg^{*(0,q)}\otimes \lieg^{p,0}=
{\mbox\rm ad}_{\Lambda_1}(q-1,1;p,0)\subseteq (q,0;p,1).
\]
Similarly, $(q,0;p,0)$ and $(q,0;p-1,1)$  are in $\ker {\mbox\rm ad}_{\Lambda_2}$, and
\begin{eqnarray*}
&&{\mbox\rm ad}_{\Lambda_2}(A)\in{\mbox\rm ad}_{\Lambda_2}\lieg^{*(0,q)}\otimes \lieg^{p,0}\\
&=&
{\mbox\rm ad}_{\Lambda_2}(q-1,1;p,0)\oplus {\mbox\rm ad}_{\Lambda_2}(q-1,1;p-1,1)
\subseteq (q,0;p+1,0)\oplus (q,0;p,1).
\end{eqnarray*}

$\adL(A)$ represents a zero class in $H^q(\lieg^{p+1})$ if and only if there exists $B$ such that
$\adL(A)=\dbar B$.  By (\ref{dbar action}), the image of $\dbar$ is contained in the components
 $(a,b;k,1)$ for some $a,b,k$. From the above observation on
 ${\mbox\rm ad}_{\Lambda_1}(A)+{\mbox\rm ad}_{\Lambda_2}(A)$, such $B$ exists only if
 $\adL(A)\in (q,0;p,1)$. Furthermore by (\ref{dbar action}) again,  $B$ is contained in the component
 $(q-1, 0; p+1, 0)$; which is contained in $\liet^{*(0, q-1)}\otimes \liet^{p+1,0}$.
Then $d_2 A$ is represented by $\adL(B)$. However, by Lemma \ref{image of Lambda},
$(q-1,0;p+1,0)$ is in the kernel of both  ${\mbox\rm ad}_{\Lambda_1}$ and
${\mbox\rm ad}_{\Lambda_2}$. Therefore,  $d_2\equiv 0$ as claimed.

\subsection{Examples}

We consider examples of holomorphic Poisson structures that satisfies the conditions in
Theorem \ref{theorem on nil}.
In particular, let $W$ be an element in $\liec^{1,0}$ and $T$ an element in $\liet^{1,0}$. Define
$\Lambda=W\wedge T$. Since $\dbar (0,0;1,1)\subseteq (1,0;0,2)$ by (\ref{dbar action}),
the dimension restriction implies that $\Lambda$ is holomorphic.
As the complex structure is abelian, $\Lambda$ is Poisson.

The first example of non-abelian nilmanifold is the Kodaira surface. As a nilmanifold,
it is covered by the product of the real three-dimensional Heisenberg group $H_3$ and
the one-dimensional trivial additive group $R^1$. Computation on this manifold was done
extensively in \cite{Poon}. It is known to have an invariant holomorphic symplectic structure.
Using the usual transformation from symplectic structure to Poisson structure, one obtains a
holomorphic Poisson structure on such a Kodaira surface such that the Poisson bi-vector field is of the type
$W\wedge T$, and hence satisfies the conditions in Theorem \ref{theorem on nil}.
We skip the details for this case, and refer readers to \cite{Poon}.

Six-dimensional 2-step nilmanifolds with abelian complex structures were studied extensively
in \cite{CFGU} \cite{MPPS-2-step} \cite{Salamon}. The underlying algebraic structure is classified.
In particular, let $R^n$ denotes the $n$-dimensional abelian group and $H_n$ the $n$-dimensional real Heisenberg
group, then the 2-step nilpotent groups with abelian complex structures are $R^6$,
$H_3\times R^3$, $H_5\times R^1$, $H_3\times H_3$, the real six-dimensional Iwasawa group $W_6$
and one additional group $P_6$ which we will describe with a little more details below.
Since $R^6$ is abelian and $H_3\times R^3$ covers the product of a Kodaira surface and an elliptic curve,
they are not new examples in a strict sense. We consider the four remaining
cases in their respective series in all admisssible dimensions, namely
$H_{2n+1}\times R^1$, $H_{2n+1}\times H_{2m+1}$, $W_{4n+2}$ and $P_{4n+2}$. We demonstrates that they will
each provide an example for Theorem \ref{theorem on nil}.

\subsubsection{$H_{2n+1}\times R^1$}
Let $\{X_{k}, Y_k, Z_1, Z_2: 1\leq k\leq n\}$ be a basis of a $(2n+2)$-dimensional vector space $\lieg$.
Define a Lie bracket by
\begin{eqnarray}\label{structure of h}
&\lbra{X_k}{Y_k}=-\lbra{Y_k}{X_k}=Z_1, \quad \lbra{X_k}{Z_1}=\lbra{Y_k}{Z_1}=0,&  \\
&\lbra{X_k}{Z_2}=\lbra{Y_k}{Z_2}=\lbra{Z_1}{Z_2}=0& \nonumber
\end{eqnarray}
for all $k$. Then $\{X_{k}, Y_k, Z_1: 1\leq k\leq n\}$  span the $(2n+1)$-dimensional Heisenberg algebra
$\lieh_{2n+1}$. With
$Z_2$, $\lieg$ is the one-dimensional trivial extension of the Heisenberg algebra. Its center is spanned by
$Z_1$ and $Z_2$. Now, define a linear map by
\[
JX_k=Y_k, \quad JY_k=-X_k, \quad JZ_1=Z_2, \quad JZ_2=-Z_1.
\]
It defines an abelian complex structure.

\subsubsection{$H_{2n+1}\times H_{2m+1}$}
Let $\{X_{k}, Y_k, Z_1: 1\leq k\leq n\}$  be a basis for the
$(2n+1)$-dimensional Heisenberg algebra $\lieh_{2n+1}$ as in (\ref{structure of h}).
Let $\{U_{\ell}, V_\ell, Z_2: 1\leq \ell\leq m\}$ be a basis for a $2m+1$-dimensional Heisenberg algebra
$\lieh_{2m+1}$ so that $\lbra{U_\ell}{V_\ell}=Z_2$ for all $\ell$.
On $\lieh_{2n+1}\oplus \lieh_{2m+1}$, its center is  spanned by $Z_1$ and $Z_2$.
Consider the linear map
\[
JX_k=Y_k, \quad JY_k=-X_k, \quad
JU_\ell=V_\ell, \quad JV_\ell=-U_\ell,
\quad
JZ_1=Z_2, \quad JZ_2=-Z_1.
\]
It defines an abelian complex structure.

\subsubsection{$W_{4n+2}$}
It is a $4n+2$-dimensional nilpotent group. On Lie algebra level, let
\begin{equation}\label{basis for w4n+2}
\{Z_1, Z_2, X_{4k+1}, X_{4k+2}, X_{4k+3}, X_{4k+4}:
0\leq k\leq n-1\}
\end{equation}
 be a basis such that the non-zero structure equations are given by
\begin{eqnarray}\label{w4n+2}
&\lbra{X_{4k+1}}{X_{4k+3}}=-\frac12 Z_1, \quad \lbra{X_{4k+1}}{X_{4k+4}}=-\frac12 Z_2, & \\
&\lbra{X_{4k+2}}{X_{4k+3}}=-\frac12 Z_2, \quad
\lbra{X_{4k+2}}{X_{4k+4}}=\frac12 Z_1. &\nonumber
\end{eqnarray}
We obtain an abelian complex structure if we insist $J^2=-$identity, and
\begin{equation}\label{j on w4n+2}
JX_{4k+1}=X_{4k+2}, \quad JX_{4k+3}=-X_{4k+4}, \quad JZ_1=-Z_2.
\end{equation}

\subsubsection{$P_{4n+2}$}
It is again a $4n+2$-dimensional nilpotent group. On Lie algebra level, it is a degeneration of
the structure in (\ref{w4n+2}). In particular, with respect to the basis
as given in (\ref{basis for w4n+2}), the non-zero structure equations are given by
\[
\lbra{X_{4k+1}}{X_{4k+2}}=-\frac12 Z_1, \quad \lbra{X_{4k+1}}{X_{4k+4}}=-\frac12 Z_2,
\quad \lbra{X_{4k+2}}{X_{4k+3}}=-\frac12 Z_2.
\]
An abelian complex structure is defined by the requirements in (\ref{j on w4n+2}).

\section{Complex Parallelizable Nilmanifolds}

A nilmanifold is complex parallelizable if there is an invariant complex structure $J$ such that for all
$X,Y$ in $\lieg$,
\begin{equation}
\lbra{X}{JY}=J\lbra{X}{Y}.
\end{equation}
In this case, $\lie g$ is the underlying real Lie algebra of a complex Lie algebra, and
\[
\lbra{\lieg^{1,0}}{\lieg^{1,0}}\subseteq {\lieg^{1,0}},
\quad \lbra{\lieg^{1,0}}{\lieg^{0,1}}=\{ 0 \}.
\]
Dually, when $\partial$ is the Chevalley-Eilenberg differential on $\lieg^{*(1,0)}$, then
 $(\oplus\lieg^{*(1,0)}, \partial)$ forms an exterior differential algebra. So is the conjugate.
Let $H^q_{\dbar}(\lieg^{*(0,1)})$ denote the $q$-th cohomology of the differential algebra
$(\oplus\lieg^{*(0,q)}, \dbar)$.
It is proved in \cite{Sakane} that the Dolbeault cohomology $H^{p,q}(M)$ is
isomorphic to the $H^q_{\dbar}(\lieg^{*(0,1)})\otimes \lieg^{*(p,0)}$.

Since the holomorphic tangent bundle for a complex parallelizable manifold is
trivialized by the complex vector fields of the algebra $\lieg^{1,0}$,
\begin{equation}
H^q(M, \Theta^p)\cong H^{0,q}(M)\otimes \lieg^{p,0}\cong H^q_{\dbar}(\lieg^{*(0,1)})\otimes\lieg^{p,0}.
\end{equation}

Suppose that $\Lambda$ is an invariant holomorphic Poisson structure on $M$.
The above observation implies that the first sheet of the Poisson double complex is given by
 \begin{equation}
 E_1^{p,q}=H^q_{\dbar}(\lieg^{*(0,1)})\otimes\lieg^{p,0}.
 \end{equation}

 Since the exterior differential of
$\lieg^{*(0,1)}$ is contained in $\lieg^{*(0,2)}$, $\adL(\oomega)=0$ for all $\oomega\in \lieg^{*(0,1)}$.
If $\sum_{j}\Oomega^j\otimes \theta_j$ represents an element in  $H^q_{\dbar}(\lieg^{*(0,1)})\otimes\lieg^{p,0}$
where $\Oomega^j$ for appropriate number of $j$ forms a basis for $H^q_{\dbar}(\lieg^{*(0,1)})$,
and $\theta_j\in \lieg^{p,0}$, then
\begin{equation}
\adL(\sum_{j}\Oomega^j\otimes \theta_j)
=\sum_{j}(-1)^q\Oomega^j\otimes \adL(\theta_j).
\end{equation}
Therefore, it represents a zero class in $H^{q+1}_{\dbar}(\lieg^{*(0,1)})\otimes\lieg^{p,0}$ if and only if
$\adL(\theta_j)=0$ for each $j$. It follows that $d_2^{p,q}\equiv 0$.
Moreover, $(\lieg^{p,0}, \adL)$ form a complex.
Let $H_{\adL}^p(\lieg^{1,0})$ be the $p$-th cohomology, then
\[
E_2^{p,q}=H^q_{\dbar}(\lieg^{*(0,1)})\otimes H_{\adL}^p(\lieg^{1,0}).
\]
Due to the degeneracy,
$H^k_{\Lambda}(M)=\oplus_{p+q=k} H^q_{\dbar}(\lieg^{*(0,1)})\otimes H_{\adL}^p(\lieg^{1,0}).$
In summary, we have
\begin{theorem}\label{parallelizable theorem}
Let $M$ be a complex parallelizable manifold with an invariant holomorphic Poisson structure $\Lambda$.
Let $\lieg^{1,0}$ the vector space of invariant holomorphic vector fields on $M$. Let
 $H^q_{\dbar}(\lieg^{*(0,1)})$ be the Chevalley-Eilenberg cohomology for the conjugation algebra
 $\lieg^{0,1}$, and $H_{\adL}^p(\lieg^{1,0})$ the cohomology of the operator $\adL=\lbra{\Lambda}{-}$
 on the exterior algebra
 of $\lieg^{1,0}$, then
 \[
 H^k_{\Lambda}(M)=\oplus_{p+q=k} H^q_{\dbar}(\lieg^{*(0,1)})\otimes H_{\adL}^p(\lieg^{1,0}).
 \]
 \end{theorem}

\subsection{Complex Three-dimensional Iwasawa Manifolds}

Let $G$ be the complex Heisenberg group of $3\times 3$ matrices of the form
\[
\mathbf{z}:
=\left(
\begin{array}{ccc}1&z_1&z_2\\
0&1&z_3\\
0&0&1
\end{array}
\right).
\]
Let $H$ be the co-compact lattice of matrices whose entries are Gaussian integers.
As a complex manifold, $G\simeq \mathbb{C}^3$.
Let $\frac{\partial}{\partial z_1}, \frac{\partial}{\partial z_2}, \frac{\partial}{\partial z_3}$ be
the corresponding coordinate global vector fields. Viewed as matrices,
\[
\frac{\partial}{\partial z_1}
=\left(
\begin{array}{ccc}0&1&0\\
0&0&0\\
0&0&0
\end{array}\right),
\qquad\frac{\partial}{\partial z_2}
=\left(
\begin{array}{ccc}0&0&1\\
0&0&0\\
0&0&0
\end{array}\right),
\qquad
\frac{\partial}{\partial z_3}
=\left(
\begin{array}{ccc}0&0&0\\
0&0&1\\
0&0&0
\end{array}
\right).
\]
Then, $\mathfrak{g}^{1,0}$ is spanned by the vector fields
\[
W_1=\frac{\partial}{\partial z_1},\qquad W_2=\frac{\partial}{\partial z_2},
\qquad W_3=\frac{\partial}{\partial z_3}+z_1\frac{\partial}{\partial z_2}
\]
which are all holomorphic. Their dual basis is given by
\[
\omega_1=dz_1,\qquad \omega_2=dz_2-z_1dz_3,\qquad \omega_3=dz_3
\]
and
$
\dbar\overline{\omega}_1=\dbar\overline{\omega}_3=0,$
$\dbar\overline{\omega}_2=-\overline{\omega}_1\wedge\overline{\omega}_3.$
Therefore,
\[
E_1^{0,0}=\langle 1\rangle,
\quad E_1^{0,1}=\langle\overline{\omega}_1,\overline{\omega}_3\rangle
\quad E_1^{0,2}=\langle\overline{\omega}_1\wedge\overline{\omega}_2,
\overline{\omega}_2\wedge\overline{\omega}_3\rangle,
\quad E_1^{0,3}=\langle\overline{\omega}_1\wedge\overline{\omega}_2\wedge\overline{\omega}_3\rangle.
\]
Also note that
\[
[W_1,W_2]=[W_2,W_3]=0,\quad [W_1,W_3]=W_2
\]
and
$
[W_i,\overline{\omega}_j]=0.$
By direct computation, the invariant holomorphic Poisson structures are given by
\[
\Lambda=aW_1\wedge W_2+ b W_2\wedge W_3,
\]
for some constants $a,b\in\mathbb{C}.$ Moreover, as we get ${\rm ad}_{\Lambda}=0$, $d_1=0$:
\begin{proposition}Let $M$ be the Iwasawa manifold,
equipped with the usual parallelizable complex structure.
If $\Lambda$ is any invariant holomorphic Poisson structure on $M$,
then the corresponding spectral sequence degenerates at the first step, and
$H^k_{\Lambda}(M)=\oplus_{p+q=k} H^q_{\dbar}(\lieg^{*(0,1)})\otimes \lieg^{1,0}$.
\end{proposition}

\section{Lichnerowicz Homomorphism}

On any complex manifold, there is a well known double complex, the Dolbeault bi-complex.
Recall that $d=\partial+\dbar$. Define
\[
B^{p,q}=C^\infty(M, \wedge^pT^{*(1,0)}M\otimes \wedge^qT^{*(0,1)}M).
\]
Due to the integrability of the complex structure, we have the double complex
\[\begin{array}
[c]{ccccc}
&  &  &  & \\
& B^{p,q+1} & \longrightarrow & B^{p+1,q+1} & \\
\overline{\partial} & \uparrow &  & \uparrow & \\
& B^{p,q} & \longrightarrow & B^{p+1,q} & \\
&  & \partial &  &
\end{array}
.
\]

Recall that $A^{p,q}=C^\infty(M, \wedge^pT^{(1,0)}M\otimes \wedge^qT^{*(0,1)}M)$.
Given a holomorphic Poisson structure $\Lambda$, contraction with $-\Lambda$ defines a map from $B^{1,0}$ to
$A^{2,0}$. We extend to the map $\phi: B^{p,q}\to A^{p,q}$ as follows.
\begin{enumerate}
\item $\phi(\omega)=-{(\Lambda\omega)}$ for all $\omega$ in $B^{1,0}$.
\item $\phi(\oomega)=\oomega$ for all $\oomega$ in $B^{0,1}$.
\item $\phi=$identity on $B^{0,0}$.
\item Extend $\phi$ to $B^{p,q}$ for $p,q\geq 1$ by $\phi (A\wedge B)=\phi(A)\wedge \phi(B)$.
\end{enumerate}

By construction, it is apparent that $\phi$ maps $B^{p,q}$ to $A^{p,q}$. The rest of this section is to prove
the following observation.

\begin{theorem}\label{homo} The map $\phi$ is a homomorphism from the Dolbeault double complex to the
Poisson double complex.
\begin{equation}\label{intertwine}
\dbar \circ \phi=\phi\circ \dbar, \quad \adL\circ \phi=\phi \circ \partial.
\end{equation}
\end{theorem}
It suffices to verify the above identities on $B^{0,0}$, $B^{1,0}$ and $B^{0,1}$.
For any function $f$, $\dbar \circ \phi f=\dbar f$. As $\dbar f$ is a $(0,1)$-form, it is clear that
$\phi \circ\dbar f=\dbar f$.
To work on the other identity, we do it locally, and there exists $(1,0)$ holomorphic vector fields
$X_j$ and $Y_j$ such that at least locally, $\Lambda=\sum_j X_j\wedge Y_j$. Then
\[
 \phi \partial f=-\sum_j (X_j\wedge Y_j)df=-\sum_j ((X_jf)Y_j-(Y_jf)X_j).
\]
On the other hand,
\[
 \adL (\phi f)= \sum_j \lbra{X_j\wedge Y_j}{f}
 =  \sum_j (X_j\wedge \lbra{Y_j}{f}-Y_j\wedge\lbra{X_j}f)
 =\sum_j((Y_jf)X_j-(X_jf)Y_j).
\]
Therefore, the identities in (\ref{intertwine}) are satisfied when they are applied to smooth functions
on the manifold.

Next, we test the identities (\ref{intertwine}) against $(0,1)$-forms, say $\oomega$.
Since the map $\phi$ is the identity map when it is restricted to the bundle of $(0,q)$-forms. Therefore,
the identity $\dbar\circ\phi(\oomega)=\phi(\dbar\oomega)$ is obviously satisfied.

\begin{lemma}\label{switching order}  For any $(2,0)$-form $\Omega$, $(1,0)$-forms $\gamma$ and $\delta$,
$\phi(\Omega)(\gamma, \delta)=\Omega(\Lambda\gamma, \Lambda\delta).$ For any $(1,1)$-form $\Gamma$,
$(1,0)$-form $\gamma$ and $(0,1)$-vector $\overline{Z}$,
$\phi(\Gamma)(\gamma, {\overline Z})=\Gamma(\Lambda\gamma, {\overline Z}).$
\end{lemma}
\bproof By linearity, it suffices to consider the case when $\Omega=\alpha\wedge\beta$ where
$\alpha$ and $\beta$ are $(1,0)$-forms. In this case, it is just a matter of definition of $\phi$ because,
\begin{eqnarray*}
&&\phi(\alpha\wedge\beta)(\gamma, \delta)=(\Lambda(\alpha)\wedge\Lambda(\beta))(\gamma, \delta)\\
&=& \Lambda(\alpha, \gamma)\Lambda(\beta, \delta)-\Lambda(\alpha, \delta)\Lambda(\beta, \gamma)
=\alpha(\Lambda\gamma)\beta(\Lambda\delta)-\alpha(\Lambda\delta)\beta(\Lambda\gamma).
\end{eqnarray*}

The proof of the second statement is similar.
\eproof

Suppose that $\gamma$ and $\overline{Z}$ are as given in the last lemma, then
\begin{eqnarray*}
&& \phi(\partial\oomega)(\gamma,{\overline{Z}})=\partial\oomega(\Lambda\gamma, {\overline{Z}})\\
&=&(\Lambda\gamma)(\oomega({\overline{Z}})-{\overline{Z}}(\oomega(\Lambda\gamma))
-\oomega(\lbra{\Lambda\gamma}{\overline{Z}})
={\cal L}_{\Lambda\gamma}\oomega({\overline{Z}})-\oomega({\cal L}_{\Lambda\gamma}{\overline{Z}})
=({\cal L}_{\Lambda\gamma}\oomega){\overline{Z}}.
\end{eqnarray*}
It is equal to $\overline{Z}(\lbra{\Lambda\gamma}{\oomega})=\adL(\oomega)(\gamma, {\overline{Z}})$
as a consequence of the observation of the next lemma.

\begin{lemma}{\rm\cite[Lemma 4, Page 10]{Brian}}\label{Brian 1}
For any $(1,0)$-vector field $X$, $(0,1)$-vector field $\overline{Z}$,
$(0,1)$-form $\oomega$ and $(1,0)$-forms  $\gamma$, and $\delta$,
\begin{eqnarray*}
\lbra{\Lambda}{X}(\gamma, \delta) &=&
-\Lambda(\omega)(\gamma(\Lambda\delta))-\gamma(\lbra{\Lambda\delta}{X})
+\delta(\lbra{\Lambda\gamma}{X}),\\
\lbra{\Lambda}{\oomega}(\gamma, \overline{Z})&=&\overline{Z}(\lbra{\Lambda\gamma}{\oomega}).
\end{eqnarray*}
\end{lemma}

Since the identities (\ref{intertwine}) holds when they are tested on $(0,1)$-forms,
we next test them on a $(1,0)$-form $\omega$. In particular,
for any $(1,0)$-forms $\omega$, $\gamma$, and $\delta$,
\begin{eqnarray}
&& \phi(\partial\omega)(\gamma, \delta)=(\partial\omega)(\Lambda\gamma, \Lambda\delta)\nonumber\\
&=&(\Lambda\gamma)(\omega(\Lambda\delta))-(\Lambda\delta)(\omega(\Lambda\gamma))
-\omega(\lbra{\Lambda\gamma}{\Lambda\delta}. \label{use Lambda}
\end{eqnarray}

Recall the following observation.
\begin{lemma}{\rm\cite[Lemma 18, Page 31]{Brian}\label{Brian 2}} \label{Poisson bracket}
Let $\Lambda$ be a bi-vector field. For any $(1,0)$-forms $\omega$, $\gamma$, and $\delta$,
\begin{eqnarray*}
&&\frac12\lbra{\Lambda}{\Lambda}(\omega, \gamma, \delta)
=\omega(\lbra{\Lambda\gamma}{\Lambda\delta})-(\Lambda\omega)(\Lambda(\gamma, \delta))\\
&+&\gamma(\lbra{\Lambda\delta}{\Lambda\omega})-(\Lambda\gamma)(\Lambda(\delta, \omega))
+\delta(\lbra{\Lambda\omega}{\Lambda\gamma})-(\Lambda\delta)(\Lambda(\omega, \gamma)).
\end{eqnarray*}
\end{lemma}

Since the bivector field $\Lambda$ is Poisson, $\lbra{\Lambda}{\Lambda}=0$.
By the last lemma, the identity (\ref{use Lambda}) is equal to
\[
-(\Lambda\omega)(\Lambda(\gamma, \delta))
+\gamma(\lbra{\Lambda\delta}{\Lambda\omega})
+\delta(\lbra{\Lambda\omega}{\Lambda\gamma}).
\]
By Lemma \ref{Brian 1}, it is equal to $-\lbra{\Lambda}{\Lambda\omega}(\gamma, \delta)$.
By definition of $\phi$, it is equal to $\adL(\phi(\omega))(\gamma, \delta)$.
As this is true for all  $\gamma$ and $\delta$,
$\adL(\phi(\omega))=\phi(\partial\omega)$ as needed.

By Lemma \ref{switching order}, for any
$(1,0)$-forms $\omega$, $\gamma$ and
$(0,1)$-vector field $\overline{Z}$,
\begin{eqnarray}
&&\phi(\dbar\omega)(\gamma, {\overline{Z}})=
(\dbar\omega)(\Lambda\gamma, {\overline{Z}})=(d\omega)(\Lambda\gamma, {\overline{Z}})
\nonumber\\
&=&(\Lambda\gamma)\omega(\overline{Z})-\overline{Z}(\omega(\Lambda\gamma))-
\omega(\lbra{\Lambda\gamma}{\overline{Z}})=-\overline{Z}(\Lambda(\gamma, \omega))-
\omega(\lbra{\Lambda\gamma}{\overline{Z}})
\nonumber\\
&=&\overline{Z}(\Lambda(\omega, \gamma))+\omega(\lbra{\overline{Z}}{\Lambda\gamma}).
\label{phidbaromega}
\end{eqnarray}

On the other hand, we treat $\dbar$ as the Lie algebroid differential for the bundle
$\overline{L}$, then
\begin{eqnarray}
&&(\dbar\phi\omega)(\gamma, \overline{Z})=-(\dbar(\Lambda\omega))(\gamma, \overline{Z})
\nonumber\\
&=&-\rho(\omega)((\Lambda\omega)(\overline{Z}))+\rho(\overline{Z})\omega(\Lambda\omega)
+\Lambda\omega(\lbra{\omega}{\overline{Z}})
\nonumber\\
&=&\overline{Z}(\omega(\Lambda\omega))
-({\cal L}_{\overline{Z}}{\omega})\Lambda\omega=\overline{Z}(\omega(\Lambda\omega))
-\overline{Z}(\omega(\Lambda\omega))+\omega(\lbra{\overline{Z}}{\Lambda\omega})
\nonumber\\
&=&\omega(\lbra{\overline{Z}}{\Lambda\omega}).
\label{dbarphiomega}
\end{eqnarray}

As a consequence of Lemma  \ref{dbar closed} and the assumption that $\Lambda$ is holomorphic, for any
$(1,0)$-forms $\omega$, $\gamma$ and
$(0,1)$-vector field $\overline{Z}$,
\begin{equation}
\overline{Z}(\Lambda(\omega,\gamma))+\omega(\lbra{\overline{Z}}{\Lambda\gamma})
=\gamma(\lbra{\overline{Z}}{\Lambda\omega}).
\end{equation}
By (\ref{phidbaromega}) and (\ref{dbarphiomega}), we conclude that the identifies in
(\ref{intertwine}) holds when it is tested against any $(1,0)$-forms.

The proof of Theorem \ref{homo} is now complete.

\subsection{Applications}

The first identity in (\ref{intertwine}) implies that the bundle map $\phi$ induces a complex linear
homomorphism
\[
\phi:H^{q}\left(  X,\Omega^{p}\right)  \longrightarrow H^{q}\left(
X,\Theta^{p}\right)  .
\]
The second identity implies that this homomorphism fits into the following
commutative diagram:
\[
\begin{array}
[c]{ccccc}
&  & d_{1}^{p,q}=ad_{\Lambda} &  & \\
& H^{q}\left(  X,\Theta^{p}\right)  & \longrightarrow & H^{q}\left(
X,\Theta^{p+1}\right)  & \\
\phi & \uparrow &  & \uparrow & \phi\\
& H^{q}\left(  X,\Omega^{p}\right)  & \longrightarrow & H^{q}\left(
X,\Omega^{p+1}\right)  & \\
&  & d_{1}^{p,q}=\partial &  &
\end{array}
,
\]
where $\Omega^p$ is the holomorphic bundle of $(p,0)$-forms.
The map is defined so that when $p=0$, $\phi$ is an identity map.

\begin{lemma}
Suppose that $X$ is holomorphic Poisson manifold such that the Fr\"{o}hlicher
spectral sequence degenerates at the $E_{1}$-level, then for spectral
sequence of  the holomorphic Poisson double complex
\[
E_{n}^{0,q}=E_{1}^{0,q}=H^{q}\left(  X,\mathcal{O}\right),
\] and $d_{n}^{0,q}\equiv 0$
for all $q$ and for all $n\geq 2$.
\end{lemma}
\bproof
As noted above, when $p=0$, the homomorphism $\phi$ is an isomorphism.
The $E_{1}^{0,q}$-term for the Fr\"{o}hlicher
spectral sequence is $H^{q}\left(  X,\mathcal{O}\right)$.
When
$d_{1}^{p,q}=\partial$ is identically equal to zero on it,
then $d_{1}^{0,q}=ad_{\Lambda}$ is
identically equal to zero on the $E_{1}^{0,q}$-term for the Poisson spectral
sequence, which is $H^{q}\left(  X,\mathcal{O}\right)$
as well. It follows that $E_{2}^{0,q}=H^{q}\left(  X,\mathcal{O}\right)  .$

Suppose that $\overline{\omega}$ is a $\overline{\partial}$-closed
$(0,q)$-form representing an element in $H^{q}\left(  X,\mathcal{O}\right)  .$
\ Since $d_{1}\overline{\omega}$ represents the zero element in $H^{q}%
(X,\Theta)$, there exists a section $A$ of $A^{1,q-1}$ such that
\[
ad_{\Lambda}\left(  \overline{\omega}\right)  =\overline{\partial}A.
\]

By definition, $d_{2}^{0,q}\left[  \overline{\omega}\right]  $ is represented
by $\adL A.$ This is an element in $A^{2,q-1}$, representing a class in
$H^{q-1}(X,\Theta^{2}).$
However, since $\phi\left(  \overline{\omega}\right)  =\overline{\omega}$, we
also have
\[
\overline{\partial}A=\adL\left(  \overline{\omega}\right)
=\adL\left(  \phi\left(  \overline{\omega}\right)  \right)
=\phi\left(  \partial\overline{\omega}\right)  =0.
\]
It follows that $\adL\left(  \overline{\omega}\right)  =\overline{\partial}A =0$.
By choosing $A=0$, we find that $d_{2}^{0,q}\left[  \overline{\omega}\right]
$ is represented by zero. Therefore, $d_{2}^{0,q}\equiv 0$ for all $q$.

For $n\geq 2$ and for all $q$, since
\[
E_{n+1}^{0,q}=\ker d_n^{0,q}:E_n^{0,q}\rightarrow E_n^{n, q-n+1}.
\]
It is now apparent that by induction, $E_{n+1}^{0,q}=E_n^{0,q}$.
After all, we have seen that if $\oomega$ is in $H^{q}\left(  X,\mathcal{O}\right)$,
then under the current assumption, $\dbar\oomega+\adL(\oomega)=0$.
\eproof

\begin{theorem}
The holomorphic Poisson spectral sequence on any compact holomorphic Poisson
surface degenerates at $E_{2}$-level.
\end{theorem}

\bproof
The first sheet of the holomorphic Poisson spectral sequence of a compact
complex surface is:
\[
\begin{array}
[c]{ccccccc}
H^{2}\left(  \mathcal{O}\right)  & \overset{d_{1}^{0,2}}{\longrightarrow} &
H^{2}\left(  \Theta\right)  & \overset{d_{1}^{1,2}}{\longrightarrow} &
H^{2}\left(  \Theta^2\right)  & \longrightarrow & 0\\
&  &  &  &  &  & \\
H^{1}\left(  \mathcal{O}\right)  & \overset{d_{1}^{0,1}}{\longrightarrow} &
H^{1}\left(  \Theta\right)  & \overset{d_{1}^{1,1}}{\longrightarrow} &
H^{1}\left(  \Theta^2\right)  & \longrightarrow & 0\\
&  &  &  &  &  & \\
H^{0}(\mathcal{O}) & \overset{d_{1}^{0,0}}{\longrightarrow} & H^{0}\left(
\Theta\right)  & \overset{d_{1}^{1,0}}{\longrightarrow} & H^{0}\left(
\Theta^2\right)  & \longrightarrow & 0
\end{array}
.
\]
Since the Fr\"{o}hlicher spectral sequence degenerates at $E_{1}$-level on any
compact complex surface, the previous lemma is applicable. Therefore, the
second sheet of the holomorphic Poisson spectral sequence is given as below.
\[
\begin{array}
[c]{ccccc}
H^{2}\left(  \mathcal{O}\right)=E_{2}^{0,2}  & E_{2}^{1,2} & E_{2}^{2,2} & & 0\\
&
\mbox{
\begin{picture}(40,20)(-20,0)
\put(-25,15){\vector(4,-1){60}}
\end{picture}
}
 &  & \\
H^{1}\left(  \mathcal{O}\right)=E_{2}^{0,1}  & E_{2}^{1,1} & E_{2}^{2,1} & & 0\\
&
\mbox{
\begin{picture}(40,20)(-20,0)
\put(-25,15){\vector(4,-1){60}}
\end{picture}
}
 &  & \\
H^{0}(\mathcal{O})=E_{2}^{0,0} & E_{2}^{1,0} & E_{2}^{2,0}  & & 0
\end{array}
.
\]
Due to a dimension restriction, the only non-trivial differential on the
second sheet are
\[
d_{2}^{0,1}:H^{1}\left(  \mathcal{O}\right)  =E_{2}^{0,1}\longrightarrow
E_{2}^{2,0}
\quad \text{ and }
\quad
d_{2}^{0,2}:H^{2}\left(  \mathcal{O}\right)
=E_{2}^{0,2}\longrightarrow E_{2}^{2,1}.
\]
By the second part of the previous lemma, these two maps are identically zero.
Therefore, $d_{2}^{p,q}\equiv0$ for all $(p,q).$
\eproof

\begin{theorem}\label{theorem through homo}
Suppose that $X$ is holomorphic Poisson manifold with complex dimension $n$.
If the Fr\"{o}hlicher spectral sequence of its complex structure degenerates
at the $E_{1}$-level, then its holomorphic Poisson spectral sequence degenerates
at the $E_{n}$-level.
\end{theorem}
\bproof
Similar to the case when $n=2$, on a complex manifold with dimension $n$, the only non-trivial terms
at the top sheet are
\[
d_{n}^{0,n-1}:E_{n}^{0,n-1}\longrightarrow
E_{n}^{n,0}\quad\text{and\quad}
d_{n}^{0,n}:=E_{n}^{0,n}\longrightarrow E_{n}^{n,1}.
\]
\eproof

\noindent{\bf Acknowledgment.} D. Grandini and Y.-S. Poon thank the Mathematical Sciences Center,
Tsinghua University, for hospitality during their visits in summer 2012.

\end{document}